\documentclass[a4paper,10pt]{article}

\usepackage[english]{babel}
\usepackage{amsmath}
\usepackage{amsthm}
\usepackage{amssymb}
\usepackage{amsfonts}
\usepackage{amsxtra}
\usepackage{color}
\usepackage[breaklinks]{hyperref}
\usepackage{enumitem}
\usepackage{tocloft}

\setenumerate{topsep=0pt, itemsep=0pt}
\setitemize{topsep=0pt, itemsep=0pt}
\setdescription{topsep=0pt, itemsep=0pt}
\setlist{noitemsep}
\renewcommand{\phi}{\varphi}
\renewcommand{\epsilon}{\varepsilon}

\theoremstyle{plain}
\newtheorem{theorem}{Theorem}[section]
\newtheorem{proposition}[theorem]{Proposition}
\newtheorem{lemma}[theorem]{Lemma}

\theoremstyle{definition}
\newtheorem{example}[theorem]{Example}
\newtheorem{remark}[theorem]{Remark}
\newtheorem*{genfct}{Generating function}
\newtheorem*{diffeq}{Differential equation}
\newtheorem*{orthnorms}{Orthonormal basis}
\newtheorem*{recrel}{Recurrence relations}
\newtheorem*{intexpr}{Integral representation}
\newtheorem*{reprodformula}{Meijer's $G$-transform}

\theoremstyle{remark}
\newenvironment{pf}{\begin{proof}[\sc Proof]}{\end{proof}}

\newcommand{\td}{\,\mathrm{d}}
\newcommand{\const}{\textup{const}}
\renewcommand{\Re}{\textup{Re }}

\numberwithin{equation}{section}

\date{}
\title{Orthogonal polynomials associated to a certain fourth order differential equation}
\author{Joachim Hilgert, Toshiyuki Kobayashi\footnote{Partially supported by Grant-in-Aid for Scientific Research (B) (18340037, 22340026), Japan Society for the Promotion of Science, and the Alexander Humboldt Foundation.}, Gen Mano, Jan M\"ollers\footnote{Partially supported by the International Research Training Group 1133 ``Geometry and Analysis of Symmetries'', and the GCOE program of the University of Tokyo.}}
\makeindex

\begin{document}

\maketitle

\begin{abstract}
We introduce orthogonal polynomials $M_j^{\mu,\ell}(x)$ as eigenfunctions of a certain self-adjoint fourth order differential operator depending on two parameters $\mu\in\mathbb{C}$ and $\ell\in\mathbb{N}_0$.

These polynomials arise as $K$-finite vectors in the $L^2$-model of the minimal unitary representations of indefinite orthogonal groups, and reduce to the classical Laguerre polynomials $L_j^\mu(x)$ for $\ell=0$.

We establish various recurrence relations and integral representations for our polynomials, as well as a closed formula for the $L^2$-norm. Further we show that they are uniquely determined as polynomial eigenfunctions.\\


\end{abstract}

\tableofcontents

\section{Introduction}

Many classical sequences of special polynomials $(P_j)_{j\in\mathbb{N}_0}$ such as Hermite polynomials, Laguerre polynomials or Jacobi polynomials have the following standard properties:
\begin{itemize}
 \item $P_j(x)$ is a polynomial of degree $j$ ($j=0,1,2,\ldots$),
 \item $P_j(x)$ is an eigenfunction of a second order differential operator.
\end{itemize}
In this article we introduce a family of polynomials $(M_j^{\mu,\ell})_{j\in\mathbb{N}_0}$ depending on two parameters $\mu\in\mathbb{C}$ and $\ell\in\mathbb{N}_0$. A distinguishing feature in our setting is that these polynomials have the following properties:
\begin{itemize}
 \item $M_j^{\mu,\ell}(x)$ is a polynomial of degree $j+\ell$ ($j=0,1,2,\ldots$),
 \item $M_j^{\mu,\ell}(x)$ is an eigenfunction of a fourth order differential operator (Theorem \ref{thm:DiffEq}).
\end{itemize}
It turns out that the polynomials $(M_j^{\mu,\ell})_{j\in\mathbb{N}_0}$ enjoy a number of good properties one typically finds for classical orthogonal polynomials: They
\begin{itemize}
 \item are unique as polynomial solutions of differential equations (Theorem \ref{thm:DiffEq}),
 \item form orthonormal bases for $L^2(\mathbb{R}_+,x^{\mu-2\ell}e^{-x}\td x)$ (Theorem \ref{thm:OrthRel}),
 \item admit integral representations (Theorem \ref{thm:IntExpr}),
 \item satisfy recurrence relations (Theorem \ref{thm:RecRel}).
\end{itemize}

Our analysis of these special polynomials $M_j^{\mu,\ell}$ is motivated by recent progress on the minimal representation of the non-compact semisimple Lie group $G=O(p,q)$ (the indefinite orthogonal group) on the Hilbert space $L^2(C)$ where $C$ is an isotropic cone in $\mathbb{R}^{p+q-2}$ (see \cite{KM07bprocj,KM07b}). In contrast to traditional analysis on homogeneous spaces, the group $G$ in our setting is too large to act geometrically on $C$. This very feature lets one expect many functional equations in the representation space arising from the action of the \lq large\rq\ group $G$ on the \lq small\rq\  representation space. In fact, some of the formulas here are predicted by unitary representation theory with $(p,q)=(\mu+3,2\ell+4)$ where $\mu$ is an odd integer.

The polynomials $M_j^{\mu,\ell}$ arise as $K$-finite vectors in the representation, and indeed generate all the $K$-types parametrized by non-negative integers $j$. In the bottom case $\ell=0$ our polynomials $M_j^{\mu,0}$ reduce to Laguerre polynomials $L_j^\mu$ and the fourth order differential operator $\mathcal{P}_{\mu,\ell}$ is of the form
\begin{equation*}
 \mathcal{P}_{\mu,0} = \mathcal{Q}_{\mu}^2 + \const
\end{equation*}
with the Laguerre operator $\mathcal{Q}_{\mu}$ having $L_j^\mu$ as eigenfunctions.

As its remarkable consequence, Laguerre polynomials (multiplied by elementary functions) give a basis of $K$-finite vectors for the minimal representation of $O(2m,4)$. Similar results were previously known only for the conformal group $O(p,2)$ (in physics terms the minimal representation of $O(4,2)$ appears as the bound states of the Hydrogen atom, and incidentally
as the quantum Kepler problem).

Our proofs of the main results rely largely on purely analytic methods, in particular the analysis of the corresponding fourth order differential equation, which we started in \cite{HKMM09a}.

The anonymous referee pointed out that in the 1940s H. L. Krall \cite{Kra40} introduced families of orthogonal polynomials associated to certain fourth order differential operators. Later A. M. Krall and L. L. Littlejohn took up the study (see e.g. \cite{LK89} for a survey). While it is not immediate that Krall's Laguerre type polynomials are in any way related to our new polynomials, this issue seems to deserve closer study. 

Notation: $\mathbb{N}_0=\{0,1,2,\ldots\}$, $\mathbb{R}_+=\{x\in\mathbb{R}:x>0\}$.

\section{Statement of the main theorems}\label{sec:thms}

We introduce the family of polynomials $\{M_j^{\mu,\ell}(x):j\in\mathbb{N}_0\}$ for $\mu\in\mathbb{C}$ and $\ell\in\mathbb{N}_0$ via meromorphic generating functions $G^{\mu,\ell}(t,x)$ of two variables $t$ and $x$ defined by
\begin{equation}
 G^{\mu,\ell}(t,x) := \frac{1}{(1-t)^{\ell+\frac{\mu+3}{2}}}\left(\frac{x}{2}\right)^{2\ell+1}
 e^{\frac{x}{2}}\widetilde{I}_{\frac{\mu}{2}}\left(\frac{tx}{2(1-t)}\right)
 \widetilde{K}_{\ell+\frac{1}{2}}\left(\frac{x}{2(1-t)}\right),\label{eq:GenFct}
\end{equation}
where $\widetilde{I}_\alpha(z)=(\frac{z}{2})^{-\alpha}I_\alpha(z)$ and $\widetilde{K}_\alpha(z)=(\frac{z}{2})^{-\alpha}K_\alpha(z)$ denote the normalized $I$- and $K$-Bessel functions. More precisely, we set
\begin{equation}
 M_j^{\mu,\ell}(x) := \frac{\Gamma(j+\mu+1)}{j!2^\mu\Gamma(j+\frac{\mu+1}{2})}\left.\frac{\partial^j}{\partial t^j}\right|_{t=0}G^{\mu,\ell}(t,x).\label{eq:DefPoly}
\end{equation}
The coefficient in front is chosen to produce suitable normalizations of the top term.

\begin{theorem}\label{thm:SpecialPolynomials}
Suppose $\mu\neq-1,-2,-3,\ldots$ and $\ell\in\mathbb{N}_0$.
\begin{enumerate}
 \item[\textup{(1)}] $M_j^{\mu,\ell}(x)$ is a polynomial of degree $j+\ell$ ($j=0,1,2,\ldots$).
 \item[\textup{(2)}] (Top term) $$M_j^{\mu,\ell}(x)=\frac{(-1)^j}{j!}x^{j+\ell}+\textup{lower order terms}.$$
 \item[\textup{(3)}] (Constant term) $$M_j^{\mu,\ell}(0)
 =\frac{2^{2\ell-\mu}\Gamma(\ell+\frac{1}{2})\Gamma(j+\mu+1)
 \left(\frac{\mu+1}{2}-\ell\right)_j}{j!\Gamma(\frac{\mu+2}{2})\Gamma(j+\frac{\mu+1}{2})},$$
 where $(a)_n=a(a+1)\cdots(a+n-1)$ is the Pochhammer symbol.
\end{enumerate}
\end{theorem}

A combinatorial formula for intermediate terms of $M_j^{\mu,\ell}(x)$ will be given in Proposition \ref{prop:ExactPoly}. Here are some further special values of the polynomials $M_j^{\mu,\ell}(x)$.

\begin{example}\label{ex:BottomCase}
\begin{enumerate}
 \item[\textup{(1)}] ($\ell=0$) The polynomials $M_j^{\mu,\ell}(x)$ for $\ell=0$ reduce to the Laguerre polynomials
  \begin{equation*}
   M_j^{\mu,0}(x) = L_j^\mu(x).
  \end{equation*}
 \item[\textup{(2)}] ($j=0$) 
The bottom of the series with $j=0$ amounts to
  \begin{equation*}
   M_0^{\mu,\ell}(x) = \sum_{k=0}^\ell{\frac{(2\ell-k)!}{k!(\ell-k)!}x^k}.
  \end{equation*}
  These polynomials appear in the explicit formula for the $K$-Bessel functions with half-integer parameter. In fact, for $\ell\in\mathbb{N}_0$ and any $\mu$ we have (see e.g. \cite[III.71~(12)]{Wat44})
  \begin{equation*}
   \widetilde{K}_{\ell+\frac{1}{2}}(z) = \sqrt{\pi}z^{-(2\ell+1)}e^{-z}M_0^{\mu,\ell}(2z).
  \end{equation*}
\end{enumerate}
\end{example}

To state the differential equation for the polynomials $M_j^{\mu,\ell}(x)$ we set $\theta:=\frac{\td}{\td x}$ and introduce the fourth order differential operator
\begin{multline*}
 \mathcal{P}_{\mu,\ell} :=
 \frac{1}{x^2}\left(\left(\theta+\mu-2\ell-1-\frac{x}{2}\right)\left(\theta+\mu-\frac{x}{2}\right)
 -\left(\frac{x}{2}\right)^2\right)\\
 \times\left(\left(\theta-2\ell-1-\frac{x}{2}\right)\left(\theta-\frac{x}{2}\right)
 -\left(\frac{x}{2}\right)^2\right)
\end{multline*}
on $\mathbb{R}_+$.

\begin{theorem}[Differential equation]\label{thm:DiffEq}
Let $\mu\neq-1,-2,-3,\ldots$ and $\ell\in\mathbb{N}_0$. For every $j\in\mathbb{N}_0$ the polynomial $M_j^{\mu,\ell}(x)$ is a solution of the fourth order differential equation
\begin{equation}
 \mathcal{P}_{\mu,\ell}u = j(j+\mu+1)u.\label{eq:DiffEqPoly}
\end{equation}
Moreover, if $\mu\geq2\ell+1$, then, up to scalar multiple, $M_j^{\mu,\ell}$ is the unique polynomial solution of this equation.
\end{theorem}

\begin{theorem}[Orthonormal basis]\label{thm:OrthRel}
For $\mu>2\ell-1$ we have
\begin{multline*}
 \int_0^\infty{M_j^{\mu,\ell}(x)M_k^{\mu,\ell}(x)x^{\mu-2\ell}e^{-x}\td x}\\
 = \left\{\begin{array}{ll}\displaystyle\frac{2\Gamma(j+\mu+1)\Gamma(j+\ell+\frac{\mu+3}{2})
 \Gamma(j-\ell+\frac{\mu+1}{2})}{j!(2j+\mu+1)\Gamma(j+\frac{\mu+1}{2})^2} & (j=k),\\0 & (j\neq k).
 \end{array}\right.
\end{multline*}
If further $\mu\geq2\ell+1$ is an odd integer, then the sequence $(M_j^{\mu,\ell})_{j\in\mathbb{N}_0}$ forms an orthogonal basis of $L^2(\mathbb{R}_+,x^{\mu-2\ell}e^{-x}\td x)$.
\end{theorem}

\begin{theorem}[Recurrence relations]\label{thm:RecRel}
Suppose $\mu\neq-1,-2,-3,\ldots$ and $\ell\in\mathbb{N}_0$. Then the polynomials $M_j^{\mu,\ell}$ are subject to the following recurrence relations:
\begin{enumerate}
\item[\textup{(1)}] The three-term recurrence relation
for $(2x\frac{\td}{\td x}-x)M_j^{\mu,\ell}(x)$:
\begin{multline*}
 (2\theta-x)M_j^{\mu,\ell}(x) = (j+1)M_{j+1}^{\mu,\ell}(x)-(\mu-2\ell+1)M_j^{\mu,\ell}(x)\\
 -(j+\mu)\frac{(2j+\mu+2\ell+1)(2j+\mu-2\ell-1)}{(2j+\mu+1)(2j+\mu-1)}M_{j-1}^{\mu,\ell}(x).
\end{multline*}
\item[\textup{(2)}] The five-term recurrence relation
for $x^2 M_j^{\mu,\ell}(x)$:
\begin{equation*}
 x^2M_j^{\mu,\ell}(x) = \sum_{k=-2}^2{a_{j,k}^{\mu,\ell}M_{j+k}^{\mu,\ell}(x)}
\end{equation*}
with coefficients
\begin{align*}
 a_{j,2}^{\mu,\ell} &= (j+1)(j+2),\\
 a_{j,1}^{\mu,\ell} &= -2(j+1)(2j+\mu+2),\\
 a_{j,0}^{\mu,\ell} &= (6j^2+6(\mu+1)j+(\mu+1)(\mu+2))\\
 & \ \ \ \ \ \ \ \ \ \ \ \ \ \ \ \ \ \ \ \ \ -\frac{4(2j^2+2(\mu+1)j+(\mu-1)(\mu+2))}{(2j+\mu-1)(2j+\mu+3)}\ell(\ell+1),\\
 a_{j,-1}^{\mu,\ell} &= -\frac{2(j+\mu)(2j+\mu)(2j+\mu+2\ell+1)(2j+\mu-2\ell-1)}{(2j+\mu-1)(2j+\mu+1)},\\
 a_{j,-2}^{\mu,\ell} &= \frac{(j+\mu-1)(j+\mu)(2j+\mu+2\ell-1)(2j+\mu-2\ell-3)}{(2j+\mu-3)(2j+\mu-1)^2}\\
 & \ \ \ \ \ \ \ \ \ \ \ \ \ \ \ \ \ \ \ \ \ \ \ \ \ \ \ \ \ \ \ \ \times\frac{(2j+\mu+2\ell+1)(2j+\mu-2\ell-1)}{(2j+\mu+1)}.
\end{align*}
\item[\textup{(3)}] The recurrence relation in $\mu$:
 \begin{multline*}
  \mu(2j+\mu-1)M_j^{\mu,\ell}(x)-2\mu(j+\mu)M_{j-1}^{\mu,\ell}(x)\\
  =(j+\mu-1)(j+\mu)M_j^{\mu-2,\ell}(x)-x^2M_{j-2}^{\mu+2,\ell}(x).
 \end{multline*}
\item[\textup{(4)}] The recurrence relation in $\ell$ ($\ell\geq1$):
 \begin{multline*}
  \nu(2j+\mu-1)M_j^{\mu,\ell}(x)-2(2\ell+1)(j+\mu)M_{j-1}^{\mu,\ell}(x)\\
  =\frac{1}{2}(2j+\mu-1)M_j^{\mu,\ell+1}(x)-(2j+\mu-1)\left(\frac{x}{2}\right)^2M_j^{\mu,\ell-1}(x).
 \end{multline*}
\item[\textup{(5)}] The recurrence relation in $\mu$ and $\ell$:
\begin{multline*}
 \left(\theta-2\ell-1-\frac{x}{2}\right)\left(2(2j+\mu-1)M_j^{\mu,\ell}(x)-4(j+\mu)M_{j-1}^{\mu,\ell}(x)\right)\\
 = 2x^2M_{j-2}^{\mu+2,\ell}(x)-(2j+\mu-1)M_j^{\mu,\ell+1}(x).
\end{multline*}
\end{enumerate}
\end{theorem}

\begin{remark}
The recurrence relation in Theorem \ref{thm:RecRel}\,(5) is useful for the actual computation of the polynomials $M_j^{\mu,\ell+1}$ for fixed $\ell$ from the polynomials $M_j^{\mu,\ell}$. In fact, since in the bottom case $\ell=0$ the polynomials are simply the Laguerre polynomials, Theorem \ref{thm:RecRel}\,(5) gives us an inductive method to calculate the series of special orthogonal polynomials. For example, for $\ell=0$ and $\mu\neq-1,-2,-3,\ldots$ we find
\begin{multline*}
 M_j^{\mu,1}(x) = \frac{2x^2}{2j+\mu-1}L_{j-2}^{\mu+2}(x)\\
 -\left(\theta-1-\frac{x}{2}\right)\left(2L_j^\mu(x)-\frac{4(j+\mu)}{2j+\mu-1}L_{j-1}^\mu(x)\right).
\end{multline*}
\end{remark}

Next we turn to integral representations of the polynomials $M_j^{\mu,\ell}$ in terms of Laguerre polynomials. Note that for these integrals it is a priory not clear that they are polynomial in $x$.

\begin{theorem}[Integral representation]\label{thm:IntExpr}
For $\Re\mu>-1$ and $x>0$ the integral
\begin{equation*}
 x^{2\ell+1}\int_0^\pi{\int_0^\infty{e^{-x(\cosh\phi-1)}
 L_j^{\ell+\frac{\mu+1}{2}}(x(\cos\theta+\cosh\phi))\sin^\mu\theta\sinh^{2\ell+1}\phi\td\phi}\td\theta}.
\end{equation*}
becomes a polynomial in $x$ of degree $j+\ell$. Further, it is equal to
\begin{equation*}
 \frac{2^\mu\ell!\Gamma(\frac{\mu+1}{2})\Gamma(j+\frac{\mu+1}{2})}{\Gamma(j+\mu+1)}M_j^{\mu,\ell}(2x).
\end{equation*}
\end{theorem}

If $\mu\geq2\ell+1$ is an odd integer, the polynomials $M_j^{\mu,\ell}$ satisfy a remarkable reproducing property with respect to Meijer's $G$-transform which is built from Meijer's $G$-function $G^{20}_{04}(z|b_1,b_2,b_3,b_4)$ (see \cite[Section 7]{HKMM09a}).

\begin{theorem}[Meijer's $G$-transform]\label{thm:ReprodFormula}
If $\mu\geq2\ell+1$ is an odd integer and $x>0$, then
\begin{multline*}
 \int_0^\infty{G^{20}_{04}\left((xy)^2\left|0,-\ell-\frac{1}{2},-\frac{\mu}{2},-\ell-\frac{\mu+1}{2}\right.
 \right)e^{-2y}y^{\mu+1}M_j^{\mu,\ell}(4y)\td y}\\
 = \frac{(-1)^j}{2}e^{-2x}x^{-(2\ell+1)}M_j^{\mu,\ell}(4x).
\end{multline*}
\end{theorem}

\section{Proofs of the main theorems}\label{sec:proofs}

The proof of Theorem \ref{thm:SpecialPolynomials} will be obtained from an explicit calculation which leads
us to a combinatorial expression of the functions $M_j^{\mu,\ell}$ in terms of Laguerre polynomials. Then the strategy is to relate the polynomials $M_j^{\mu,\ell}$ to a series of functions which was already studied thoroughly in \cite{HKMM09a}. This will help to prove Theorem \ref{thm:DiffEq} and reduce the remaining theorems to results from \cite{HKMM09a}.

We start with a proof of the fact that the functions $M_j^{\mu,\ell}$ are indeed polynomials.

\begin{proposition}\label{prop:ExactPoly}
Suppose $\mu\neq-1,-2,-3,\ldots$ and $\ell\in\mathbb{N}_0$. Then $M_j^{\mu,\ell}(x)$ ($j=0,1,2,\ldots$) is given by
\begin{align}
 M_j^{\mu,\ell}(x) &= \frac{\Gamma(j+\mu+1)}{\Gamma(j+\frac{\mu+1}{2})}\sum_{k=0}^j{\sum_{i=0}^{\ell-k}{(-1)^k\frac{\Gamma(j-k+\frac{\mu+1}{2})(2\ell-i)!} {k!\Gamma(j-k+\mu+1)(\ell-i-k)!i!}L_{j-k}^\mu(x)x^i}}\notag\\
 &= \sum_{k=0}^{j+\ell}{\beta^{\mu,\ell}_{j,k}x^k},\label{eq:ExplPoly}
\end{align}
where $L_n^\alpha$ denotes the Laguerre polynomials and
\begin{multline*}
 \beta^{\mu,\ell}_{j,k} = \frac{\Gamma(j+\mu+1)}{\Gamma(j+\frac{\mu+1}{2})}\sum_{(m,n)\in S_{j,k}^{\mu,\ell}}{(-1)^{m+n}\frac{\Gamma(j-m+\frac{\mu+1}{2})}{\Gamma(n+\mu+1)}}\\
 \times{\frac{(2\ell+n-k)!}{m!n!(k-n)!(j-m-n)!(\ell+n-k-m)!}}
\end{multline*}
with
\begin{equation}
 S_{j,k}^{\mu,\ell} = \left\{(m,n)\in\mathbb{N}_0^2:\begin{array}{c}0\leq n\leq j-m\\0\leq k-n\leq\ell-m\end{array}\right\}.
 \label{eq:Sjk}
\end{equation}
\end{proposition}

\begin{pf}
Let us first assume $\Re\mu>-1$. The $K$-Bessel functions with half-integer parameter can be written explicitly as (see e.g. \cite[III.71~(12)]{Wat44})
\begin{equation}
 \widetilde{K}_{\ell+\frac{1}{2}}(z) = \sqrt{\pi}z^{-(2\ell+1)}e^{-z}\sum_{i=0}^{\ell}{\frac{(2\ell-i)!}{(\ell-i)!\cdot i!}(2z)^i}.
 \label{eq:ExplKBessel}
\end{equation}
Using the following integral representation for the $I$-Bessel function (cf. \cite[III.71~(9)]{Wat44})
\begin{align*}
 \widetilde{I}_\alpha(z) &= \frac{1}{\sqrt{\pi}\Gamma(\alpha+\frac{1}{2})}\int_0^\pi{e^{-x\cos\phi}\sin^{2\alpha}\phi\td\phi}, & \Re\alpha>-\frac{1}{2},
\end{align*}
we obtain
\begin{align*}
 G^{\mu,\ell}(t,x) ={}& \frac{1}{(1-t)^{\ell+\frac{\mu+3}{2}}}\left(\frac{x}{2}\right)^{2\ell+1}e^{\frac{x}{2}}\widetilde{I}_{\frac{\mu}{2}}\left(\frac{tx}{2(1-t)}\right)\widetilde{K}_{\ell+\frac{1}{2}}\left(\frac{x}{2(1-t)}\right)\\
 ={}& \frac{1}{\Gamma(\frac{\mu+1}{2})} \int_0^\pi{\frac{1}{(1-t)^{\frac{\mu-1}{2}+1}}e^{-\frac{tx}{2(1-t)}(\cos\phi+1)}\sin^\mu\phi\td\phi}\\
 & \ \ \ \ \ \ \ \ \ \ \ \ \ \ \ \ \ \ \ \ \ \ \ \ \ \ \ \ \ \ \ \ \ \ \ \ \ \ \ \ \times\sum_{i=0}^{\ell}{\frac{(2\ell-i)!}{(\ell-i)!\cdot i!}x^i(1-t)^{\ell-i}}.
\end{align*}
Next, we compute the derivatives of the first factor with respect to $t$ at $t=0$. Using the formula of the generating function for the Laguerre polynomials (see e.g. \cite[(6.2.4)]{AAR99})
\begin{equation}
 \sum_{n=0}^\infty{L_n^\alpha(x)t^n} = \frac{1}{(1-t)^{\alpha+1}}e^{-\frac{tx}{1-t}}\label{eq:GenFctLag}
\end{equation}
we find that
\begin{align*}
 & \left.\frac{\partial^j}{\partial t^j}\right|_{t=0}\left[\int_0^\pi{\frac{1}{(1-t)^{\frac{\mu-1}{2}+1}}e^{-\frac{tx}{2(1-t)}(\cos\phi+1)}\sin^\mu\phi\td\phi}\right]\\
 ={}& j!\int_0^\pi{L_j^{\frac{\mu-1}{2}}\left(\frac{x}{2}(\cos\phi+1)\right)\sin^\mu\phi\td\phi}\\
 \intertext{and substituting $y=\frac{1}{2}(\cos\phi+1)$ yields}
 ={}& j!2^\mu\int_0^1{(1-y)^{\frac{\mu-1}{2}}y^{\frac{\mu-1}{2}}L_j^{\frac{\mu-1}{2}}(xy)dy}\\
 ={}& \frac{j!2^\mu\Gamma(j+\frac{\mu+1}{2})\Gamma(\frac{\mu+1}{2})}{\Gamma(j+\mu+1)}L_j^\mu(x),
\end{align*}
where the last equality is the integral formula \cite[16.6~(5)]{EMOT54a}. Now we can compute the Taylor coefficients of $G^{\mu,\ell}(t,x)$ at $t=0$ explicitly as follows
\begin{align*}
 & \left.\frac{\partial^j}{\partial t^j}\right|_{t=0}G^{\mu,\ell}(t,x)\\
 ={}& \frac{1}{\Gamma(\frac{\mu+1}{2})} \sum_{k=0}^j{{j\choose k}\left.\frac{\partial^{j-k}}{\partial t^{j-k}}\right|_{t=0}\left[\int_0^\pi{\frac{1}{(1-t)^{\frac{\mu-1}{2}+1}}e^{-\frac{tx}{2(1-t)}(\cos\phi+1)}\sin^\mu\phi\td\phi}\right]}\\
 & \ \ \ \ \ \ \ \ \ \ \ \ \ \ \ \ \ \ \ \ \ \ \ \ \ \ \ \ \ \ \ \ \ \ \ \ \ \ \ \ \ \ \ \ \times\left.\frac{\partial^k}{\partial t^k}\right|_{t=0}\left[\sum_{i=0}^{\ell}{\frac{(2\ell-i)!}{(\ell-i)!\cdot i!}x^i(1-t)^{\ell-i}}\right]\\
 ={}& \sum_{k=0}^j{\frac{j!2^\mu\Gamma(j-k+\frac{\mu+1}{2})}{k!\Gamma(j-k+\mu+1)}L_{j-k}^\mu(x)}\\
 & \ \ \ \ \ \ \ \ \ \ \ \ \ \ \ \ \ \ \ \ \ \ \ \ \ \ \ \ \ \ \times\sum_{i=0}^{\ell}{\frac{(2\ell-i)!}{(\ell-i)!\cdot i!}x^i(-1)^k\left(\ell-i\right)\cdots\left(\ell-i-k+1\right)}\\
 ={}& \sum_{k=0}^j{\sum_{i=0}^{\ell-k}{(-1)^k\frac{j!2^\mu\Gamma(j-k+\frac{\mu+1}{2})(2\ell-i)!}{k!\Gamma(j-k+\mu+1)(\ell-i-k)!i!}L_{j-k}^\mu(x)x^i}}.
\end{align*}
This gives the first expression for $M_j^{\mu,\ell}(x)$. Inserting the explicit formula (cf. \cite[(6.2.2)]{AAR99})
\begin{equation*}
 L_n^\alpha(x) = \frac{\Gamma(n+\alpha+1)}{n!}\sum_{k=0}^n{(-1)^k{n\choose k}\frac{x^k}{\Gamma(k+\alpha+1)}}
\end{equation*}
one obtains the expressions for the coefficients $\beta^{\mu,\ell}_{j,k}$ given in the proposition. Since these clearly have meromorphic continuation for $\mu\in\mathbb{C}$ with poles at most at $\mu=-1,-2,-3,\ldots$, the claim follows.
\end{pf}

\begin{proof}[\sc Proof of Theorem \ref{thm:SpecialPolynomials}]
It remains to compute top and bottom term of the polynomials $M_j^{\mu,\ell}(x)$. For $k=j+\ell$ the set $S_{j,k}^{\mu,\ell}$ defined in \eqref{eq:Sjk} only contains the tuple $(0,j)$ and we
obtain the top term
\begin{equation*}
 \beta^{\mu,\ell}_{j,j+\ell}=\frac{(-1)^j}{j!}.
\end{equation*}
To calculate the bottom term $M_j^{\mu,\ell}(0)$ simply observe that (see e.g. \cite[Chapter III.7]{Wat44})
\begin{equation*}
 \widetilde{I}_{\frac{\mu}{2}}(0) = \frac{1}{\Gamma(\frac{\mu+2}{2})}\quad \text{ and }\quad\left.\left(\frac{x}{2}\right)^{2\ell+1}\widetilde{K}_{\ell+\frac{1}{2}}(x)\right|_{x=0} = \frac{\Gamma(\ell+\frac{1}{2})}{2},
\end{equation*}
so that
\begin{equation*}
 G^{\mu,\ell}(t,0)
 = \frac{2^{2\ell+1}\Gamma(\ell+\frac{1}{2})}{2\Gamma(\frac{\mu+2}{2})} (1-t)^{\ell-\frac{\mu+1}{2}}
 = \sum_{j=0}^\infty{\frac{2^{2\ell}\Gamma(\ell+\frac{1}{2})(\frac{\mu+1}{2}-\ell)_j}{j!\Gamma(\frac{\mu+2}{2})}t^j}.
\end{equation*}
Together with \eqref{eq:DefPoly} this proves the claim.
\end{proof}

In order to be able to apply results from \cite{HKMM09a} in our context we need two observations relating the polynomials $M_j^{\mu,\ell}(x)$ and the differential operators $\mathcal{P}_{\mu,\ell}$ to corresponding objects introduced in \cite{HKMM09a}.

\begin{lemma}\label{lem:RelToLambda}
In the notation of \cite{HKMM09a} we have
\begin{equation*}
 \Lambda_{2,j}^{\mu,2\ell+1}(x) = \frac{2^\mu\Gamma(j+\frac{\mu+1}{2})}{\Gamma(j+\mu+1)}x^{-(2\ell+1)}e^{-x}M_j^{\mu,\ell}(2x).
\end{equation*}
\end{lemma}

\begin{pf}
Directly from \cite[equations (3.2) and (4.2)]{HKMM09a}.
\end{pf}

Next, we relate the operator $\mathcal{P}_{\mu,\ell}$ to the fourth order differential operator
\begin{equation}\label{eq:Dmunu}
 \mathcal{D}_{\mu,\nu} = \frac{1}{x^2}\left((\theta+\nu)(\theta+\mu+\nu)-x^2\right)\left(\theta(\theta+\mu)-x^2\right)-\frac{(\mu-\nu)(\mu+\nu+2)}{2}
\end{equation}
introduced in \cite{HKMM09a}.

\begin{lemma}\label{lem:RelToDmunu}
For $f\in\mathcal{C}^\infty(\mathbb{R}_+)$ we have
\begin{multline*}
 \mathcal{D}_{\mu,2\ell+1}(x^{-(2\ell+1)}e^{-x}f(2x))\\
 = x^{-(2\ell+1)}e^{-x}\left(\left(4\mathcal{P}_{\mu,\ell} +\frac{(\mu-2\ell-1)(\mu+2\ell+3)}{2}\right)f\right)(2x).
\end{multline*}
\end{lemma}

\begin{pf}
Since $\mathcal{D}_{\mu,\nu}$ is symmetric in $\mu$ and $\nu$ by \cite[Proposition 2.1~(1)]{HKMM09a}, we can rewrite the expression \eqref{eq:Dmunu} for $\mathcal{D}_{\mu,\nu}$ as
\begin{equation*}
 \mathcal{D}_{\mu,\nu} = \frac{1}{x^2}\left((\theta+\mu)(\theta+\mu+\nu)-x^2\right)\left(\theta(\theta+\nu)-x^2\right) +\frac{(\mu-\nu)(\mu+\nu+2)}{2}.
\end{equation*}
Then for $\nu=2\ell+1$ the claim follows by using the commutator relation
\begin{equation*}
 \left[\theta,x^{-\nu}e^{-x}\right] = -(\nu+x)x^{-\nu}e^{-x}.\qedhere
\end{equation*}
\end{pf}

Now, in view of Lemmas \ref{lem:RelToLambda} and \ref{lem:RelToDmunu}, it follows from \cite[Theorem 4.6]{HKMM09a} that $M_j^{\mu,\ell}(x)$ solves the differential equation \eqref{eq:DiffEqPoly}. The crucial observation for the proof of uniqueness in Theorem \ref{thm:DiffEq} is that the subspace of solutions of \eqref{eq:DiffEqPoly} bounded near $x=0$ can be wholly described in terms of the polynomials $M_j^{\mu,\ell}(x)$. Note, that the operator $\mathcal{P}_{\mu,\ell}$ is invariant under the transformation $f(x)\mapsto\widetilde{f}(x)=e^xf(-x)$, i.e.
\begin{equation*}
 \mathcal{P}_{\mu,\ell}\widetilde{f} = \widetilde{\mathcal{P}_{\mu,\ell}f}.
\end{equation*}
Then clearly $N_j^{\mu,\ell}(x):=\widetilde{M_j^{\mu,\ell}}(x)=e^xM_j^{\mu,\ell}(-x)$ is another non-trivial solution of \eqref{eq:DiffEqPoly}. Since $N_j^{\mu,\ell}(x)$ grows exponentially as $x\rightarrow\infty$ and $M_j^{\mu,\ell}(x)$ just polynomially, the two solutions are linearly independent. We even have the following lemma:

\begin{lemma}\label{lem:FundSys}
Suppose $\mu\geq2\ell+1$. Then the subspace of solutions of \eqref{eq:DiffEqPoly} which are bounded near $x=0$ is two-dimensional and spanned by the functions $M_j^{\mu,\ell}(x)$ and $N_j^{\mu,\ell}(x)$.
\end{lemma}

\begin{pf}
By the previous considerations it remains to show that the subspace of solutions of \eqref{eq:DiffEqPoly} which are bounded near $x=0$ is two-dimensional.

We note first that the differential operator $\mathcal{P}_{\mu,\ell}$ has a regular singularity at $x=0$ with characteristic exponents $\{0,-\mu,2\ell+1,2\ell+1-\mu\}$. In fact, an easy computation shows that
\begin{equation*}
 x^2\mathcal{P}_{\mu,\ell}\equiv(\theta+\mu-2\ell-1)(\theta+\mu)(\theta-2\ell-1)\theta\ \ \ \ \ (\textup{mod }x\cdot\mathbb{C}[x,\theta]),
\end{equation*}
where $\mathbb{C}[x,\theta]$ denotes the left $\mathbb{C}[x]$-module generated by $1,\theta,\theta^2,\ldots$ in the Weyl algebra $\mathbb{C}[x,\frac{\td}{\td x}]$. Therefore, the differential equation $\mathcal{P}_{\mu,\ell}u=\lambda u$ is of regular singularity at $x=0$, and its characteristic equation is given by
\begin{equation*}
 (s+\mu-2\ell-1)(s+\mu)(s-2\ell-1)s=0.
\end{equation*}

Since $2\ell+1>0\geq2\ell+1-\mu>-\mu$ for $\mu\geq2\ell+1$, the theory of regular singularities (see e.g. \cite[Chapter 4]{CL55}) assures that the subspace of solutions of \eqref{eq:DiffEqPoly} bounded near $x=0$ is two-dimensional. In fact, for $\mu>2\ell+1$ there are exactly two non-negative exponents, namely $2\ell+1$ and $0$. Thus, the subspace of solutions which are bounded near $x=0$ is spanned by functions $f_1(x)$ and $f_2(x)$ with $f_1(x)\sim x^{2\ell+1}$ and $f_2(x)\sim x^0=1$ for $x\rightarrow0$. If $\mu=2\ell+1$, then there are three non-negative exponents, but the exponent $0$ has multiplicity two. Hence, the subspace of solutions with asymptotic behavior at $x=0$ given by these exponents is three-dimensional and spanned by functions $f_1(x)$, $f_2(x)$ and $f_3(x)$ with asymptotic behavior at $x=0$ given by $f_1(x)\sim x^{2\ell+1}$, $f_2(x)\sim x^0=1$ and $f_3(x)\sim\log(x)$. But $f_3(x)$ is clearly not bounded near $x=0$, so again the subspace of solutions bounded near $x=0$ is two-dimensional.
\end{pf}

\begin{proof}[\sc Proof of Theorem \ref{thm:DiffEq}]
Only the uniqueness remains to be shown. To do that, observe that by Lemma \ref{lem:FundSys} every polynomial solution of \eqref{eq:DiffEqPoly} has to be a linear combination of $M_j^{\mu,\ell}(x)$ and $N_j^{\mu,\ell}(x)$. But $N_j^{\mu,\ell}$ grows exponentially as $x\rightarrow\infty$ and hence cannot be a polynomial. This leaves $M_j^{\mu,\ell}$ as the only polynomial solution of \eqref{eq:DiffEqPoly} (up to scalar multiples).
\end{proof}

With these preparations the remaining results of Section \ref{sec:thms} are immediate consequences of results from \cite{HKMM09a}:

\begin{proof}[\sc Proof of Theorem \ref{thm:OrthRel}]
This follows from \cite[Corollaries 4.8 and 6.2]{HKMM09a}.
\end{proof}

\begin{proof}[\sc Proof of Theorem \ref{thm:RecRel}]
In view of Lemma \ref{lem:RelToLambda} the recurrence relations
follow from \cite[Propositions 6.1, 6.4 and 6.6]{HKMM09a}.
\end{proof}

\begin{proof}[\sc Proof of Theorem \ref{thm:IntExpr}]
This is essentially a reformulation of \cite[Theorem 5.1~(1)]{HKMM09a}.
\end{proof}

\begin{proof}[\sc Proof of Theorem \ref{thm:ReprodFormula}]
The stated formula follows immediately from \cite[Theorem 7.2]{HKMM09a}.
\end{proof}

\section{Comparison with Laguerre polynomials}

We saw in Example \ref{ex:BottomCase} that the polynomials $M_j^{\mu,\ell}$ are Laguerre polynomials if $\ell=0$. Thus our results on the $M_j^{\mu,\ell}$ specialize to results on Laguerre polynomials. In this section we examine how these results are related to standard results on Laguerre polynomials.

We start by collecting a number of formulas for Laguerre polynomials (cf. \cite[Chapter 6.2]{AAR99} and \cite[\textbf{II}, Chapter 10.12]{EMOT81}). The polynomial $L_n^\alpha(x)$ is defined by
\begin{equation}
 L_n^\alpha(x) = \sum_{k=0}^n{\frac{(-1)^k}{k!}{n+\alpha\choose n-k}x^k}.\label{eq:DefLag}
\end{equation}
We start with the three independent formulas
\begin{align}
 xL_n^\alpha &= -(n+1)L_{n+1}^\alpha+(2n+\alpha+1)L_n^\alpha-(n+\alpha)L_{n-1}^\alpha,\label{eq:LaguerreId1}\\
 \frac{\td}{\td x}L_n^\alpha &= \frac{\td}{\td x}L_{n-1}^\alpha-L_{n-1}^\alpha,\label{eq:LaguerreId4}\\
 \frac{\td}{\td x}L_n^\alpha &= -L_{n-1}^{\alpha+1}.\label{eq:LaguerreId3}
\end{align}
These imply the four additional identities
\begin{align}
 \theta L_n^\alpha &= nL_n^\alpha-(n+\alpha)L_{n-1}^\alpha,\label{eq:LaguerreId2}\\
 \left(x\frac{\td^2}{\td x^2}+(\alpha+1-x)\frac{\td}{\td x}+n\right)L_n^\alpha &= 0,\label{eq:LagDiffEq2}\\
 xL_n^{\alpha+1} &= (n+\alpha+1)L_n^\alpha-(n+1)L_{n+1}^\alpha,\label{eq:LaguerreId5}\\
 L_{n-1}^\alpha &= L_n^\alpha-L_n^{\alpha-1}.\label{eq:LaguerreId6}
\end{align}
We will also need the summation formula
\begin{equation}
 \sum_{j=0}^n{\frac{\td^j}{\td x^j}L_n^\alpha(x)} = \sum_{k=0}^n{(-1)^jL_{n-j}^{\alpha+j}(x)} = L_n^{\alpha-1}(x).\label{eq:LagSumm}
\end{equation}
which follows from \eqref{eq:LaguerreId3} and \eqref{eq:DefLag} by a simple calculation.

Now we can examine how the various types of results from Section \ref{sec:thms} specialize to Laguerre polynomials.

\begin{genfct}
For $\ell=0$ the expression \eqref{eq:GenFct} for the generating function $G^{\mu,\ell}(t,x)$ of the polynomials $M_j^{\mu,\ell}(x)$ can be simplified using (see \cite[3.71~(13)]{Wat44})
\begin{equation*}
 \widetilde{K}_{\frac{1}{2}}(x) = \frac{\sqrt{\pi}}{x}e^{-x}.
\end{equation*}
Then for $\ell=0$ combining the equations \eqref{eq:GenFct} and \eqref{eq:DefPoly} yields
\begin{equation}
 \sum_{j=0}^\infty{\frac{\Gamma(j+\frac{\mu+1}{2})}{\Gamma(j+\mu+1)}t^jL_j^\mu(x)} = \frac{\sqrt{\pi}}{2^{\mu}(1-t)^{\frac{\mu+1}{2}}}e^{-\frac{tx}{2(1-t)}}\widetilde{I}_{\frac{\mu}{2}}\left(\frac{tx}{2(1-t)}\right)\label{eq:LaguerreGenFct2}
\end{equation}
To see that \eqref{eq:LaguerreGenFct2} is in agreement with the standard generating function for the Laguerre polynomials given by \eqref{eq:GenFctLag} one can use the following two formulas (\cite[\textbf{I}, 6.12~(5)]{EMOT81} and \cite[\textbf{II}, 7.2.2~(12)]{EMOT81})
\begin{align*}
 \sum_{n=0}^\infty{\frac{(c)_n}{(\alpha+1)_n}t^nL_n^\alpha(x)} &= (1-t)^{-c}{_1F_1}\left(c;\alpha+1;-\frac{tx}{1-t}\right),\\
 \widetilde{I}_\alpha(x) &= \frac{e^{-x}}{\Gamma(\alpha+1)}{_1F_1}\left(\alpha+\frac{1}{2};2\alpha+1;2x\right),
\end{align*}
where ${_1F_1}(a;b;x)$ denotes the hypergeometric function.
\end{genfct}

\begin{diffeq}
For $\ell=0$ the fourth order differential operator $\mathcal{P}_{\mu,\ell}$ essentially degenerates to the square of a second order operator. In fact,
\begin{align*}
 \mathcal{P}_{\mu,0} &= \mathcal{Q}_\mu^2-\left(\frac{\mu+1}{2}\right)^2\\
 \intertext{with}
 \mathcal{Q}_\mu &= \frac{1}{x}\left(\theta^2+(\mu-x)\theta-\frac{\mu+1}{2}x\right) = x\frac{\td^2}{\td x^2}+(\mu+1-x)\frac{\td}{\td x}-\frac{\mu+1}{2}.
\end{align*}
But the Laguerre differential equation \eqref{eq:LagDiffEq2} shows that the polynomials $M_j^{\mu,0}=L_j^\mu$ are eigenfunctions of $\mathcal{Q}_\mu$ for the eigenvalue $-(j+\frac{\mu+1}{2})$. The fourth order differential equation \eqref{eq:DiffEqPoly} follows from this by applying $\mathcal{Q}_\mu$ twice.
\end{diffeq}

\begin{orthnorms}
In the special case $\ell=0$, we have $M_j^{\mu,\ell}(x)=L_j^\mu(x)$ 
(see Example~\ref{ex:BottomCase}~(1)), and therefore
 Theorem \ref{thm:OrthRel} reduces to the well-known fact that the Laguerre polynomials $(L_j^\mu)_{j\in\mathbb{N}_0}$ form a complete orthogonal system of $L^2(\mathbb{R}_+,x^\mu e^{-x}\td x)$ with norms
\begin{equation}
 \|L_j^\mu\|_{L^2(\mathbb{R}_+,x^\mu e^{-x}\td x)}^2 = \frac{\Gamma(j+\mu+1)}{j!}.\label{eq:norm-formula}
\end{equation}
Our result in this special case is slightly weaker
in the sense that we assumed $\mu\in2\mathbb{Z}+1$ 
for the completeness of the sequence $(M_j^{\mu,\ell})_{j\in\mathbb{N}_0}$
in $L^2(\mathbb{R}_+,x^{\mu-2\ell}e^{-x}\td x)$ whereas this remains true for the Laguerre polynomials
for arbitrary $\mu>-1$ (see e.g. \cite[Chapter 6]{AAR99}).
\end{orthnorms}

\begin{recrel}
We examine three of the different recurrence relations given in Theorem \ref{thm:RecRel}.
\begin{enumerate}
 \item[\textup{(1)}] For $\ell=0$ the three-term recurrence relation of Theorem \ref{thm:RecRel}~(1) simplifies to
  \begin{equation*}
   (2\theta-x)L_j^\mu=(j+1)L_{j+1}^\mu-(\mu+1)L_j^\mu-(j+\mu)L_{j-1}^\mu.
  \end{equation*}
  One can use the identities \eqref{eq:LaguerreId1} and \eqref{eq:LaguerreId2} to give an independent proof of this equation.
 \item[\textup{(2)}] In the case that $\ell=0$ the five-term recurrence relation of Theorem \ref{thm:RecRel}~(2) arises from the three-term recurrence relation
  \eqref{eq:LaguerreId1}. The existence of this three term recurrence relation is predicted by the general theory of orthogonal polynomials. In fact, every sequence of orthogonal polynomials starting in degree $0$ and increasing degree by $1$ in every step satisfies three-term recurrence relations for the multiplication by $x$ (see e.g. \cite[Chapter 5.2]{AAR99}). The polynomial $M_j^{\mu,\ell}$, however, is of degree $j+\ell$ and hence the degrees of the polynomials in this sequence start with $0$ only in the case where $\ell=0$. Already for $\ell=1$ it is possible to prove that there are no three-term recurrence relations for the multiplication by $x$. Indeed, if a three-term recurrence relation existed, there would be constants $a,b\in\mathbb{C}$ such that
  \begin{equation}\label{eq:3TermRecRelEx}
   xM_0^{\mu,1}(x) = aM_1^{\mu,1}(x) + bM_0^{\mu,1}(x).
  \end{equation}
  Since in this case the polynomials have the specific form
  \begin{align*}
   M_0^{\mu,1}(x) &= x+2,\\
   M_1^{\mu,1}(x) &= -x^2+(\mu-1)x+2(\mu-1),
  \end{align*}
  the equation \eqref{eq:3TermRecRelEx} is equivalent to
  \begin{equation*}
   x(x+2) = a(-x^2+(\mu-1)x+2(\mu-1))+b(x+2)
  \end{equation*}
  which cannot hold for all $x$.
 \item[\textup{(3)}] The formula in Theorem \ref{thm:RecRel}~(3) for $\ell=0$ reduces to
  \begin{equation*}
   x^2L_{j-2}^{\mu+2} = (j+\mu-1)(j+\mu)L_j^{\mu-2} - \mu(2j+\mu-1)L_j^\mu +2\mu(j+\mu)L_{j-1}^\mu.
  \end{equation*}
  This identity also follows by applying \eqref{eq:LaguerreId5} twice and \eqref{eq:LaguerreId6} three times (in this order):
  \begin{align*}
   x^2L_{j-2}^{\mu+2} ={}& x\left[(j+\mu)L_{j-2}^{\mu+1}-(j-1)L_{j-1}^{\mu+1}\right]\\
   ={}& (j+\mu-1)(j+\mu)L_{j-2}^\mu-2(j-1)(j+\mu)L_{j-1}^\mu+j(j-1)L_j^\mu\\
   ={}& -(j+\mu-1)(j+\mu)L_{j-1}^{\mu-1}+(-j+\mu+1)(j+\mu)L_{j-1}^\mu+j(j-1)L_j^\mu\\
   ={}& (j+\mu-1)(j+\mu)L_j^{\mu-2}-(j+\mu-1)(j+\mu)L_j^{\mu-1}\\
   &\ \ \ \ \ \ \ \ \ \ \ \ \ \ \ \ \ \ \ \ \ \ \ \ \ \ \ \ \ \ \ \ \ \ +(-j+\mu+1)(j+\mu)L_{j-1}^\mu+j(j-1)L_j^\mu\\
   ={}& (j+\mu-1)(j+\mu)L_j^{\mu-2} - \mu(2j+\mu-1)L_j^\mu +2\mu(j+\mu)L_{j-1}^\mu.
  \end{align*}
\end{enumerate}
\end{recrel}

\begin{intexpr}
In the case where $\ell=0$ the integral representation of Theorem \ref{thm:IntExpr} amounts to
\begin{multline*}
 L_j^\mu(x) = \frac{\Gamma(j+\mu+1)}{2^\mu\Gamma(\frac{\mu+1}{2})\Gamma(j+\frac{\mu+1}{2})}\\
 \int_0^\pi{\left(\int_0^\infty{e^{-x(\cosh\phi-1)} L_j^{\frac{\mu+1}{2}}(x(\cos\theta+\cosh\phi))x\sinh\phi\td\phi}\right)\sin^\mu\theta\td\theta}.
\end{multline*}
Iterated integration by parts for the inner integral, using
\begin{equation*}
 \frac{\partial}{\partial\phi}e^{-x(\cosh\phi-1)}=-x e^{-x(\cosh\phi-1)}\sinh\phi,
\end{equation*}
yields
\begin{multline*}
 L_j^\mu(x) = \frac{\Gamma(j+\mu+1)}{2^\mu\Gamma(\frac{\mu+1}{2})\Gamma(j+\frac{\mu+1}{2})}\\
 \int_0^\pi{\left[L_j^{\frac{\mu+1}{2}}+\cdots+\frac{\td^j}{\td x^j}L_j^{\frac{\mu+1}{2}}\right](x(\cos\theta+1))\sin^\mu\theta\td\theta}.
\end{multline*}
Using the summation formula \eqref{eq:LagSumm} and the substitution $y:=\cos\theta+1$ we finally obtain the integral formula \cite[16.6~(5)]{EMOT54a}
\begin{align*}
 L_j^\mu(x)
 &= \frac{\Gamma(j+\mu+1)}{2^\mu\Gamma(\frac{\mu+1}{2})\Gamma(j+\frac{\mu+1}{2})}\int_0^\pi{L_j^{\frac{\mu-1}{2}}(x(\cos\theta+1))\sin^\mu\theta\td\theta}\\
 &= \frac{\Gamma(j+\mu+1)}{2^\mu\Gamma(\frac{\mu+1}{2})\Gamma(j+\frac{\mu+1}{2})}\int_0^1{L_j^{\frac{\mu-1}{2}}(xy)y^{\frac{\mu-1}{2}}(1-y)^{\frac{\mu-1}{2}}\td y}.
\end{align*}
\end{intexpr}

\begin{reprodformula}
For $\ell=0$ the $G$-function appearing in Theorem \ref{thm:ReprodFormula} reduces to a $J$-Bessel function, namely (cf. \cite[\textbf{I}, 5.6~(11)]{EMOT81})
\begin{align*}
 G^{20}_{04}\left(t\left|0,-\frac{1}{2},-\frac{\mu}{2},-\frac{\mu+1}{2}\right.\right) &= t^{-\frac{\mu}{4}}J_\mu(4t^{\frac{1}{4}}), & t>0.
\end{align*}
Then (after a suitable substitution) the reproducing property of Theorem \ref{thm:ReprodFormula} can be rewritten as
\begin{equation*}
 \int_0^\infty{J_\mu(xy)(xy)^{\frac{1}{2}}y^{\mu+\frac{1}{2}}e^{-\frac{1}{2}y^2}L_j^\mu(y^2)\td y} = (-1)^jx^{\mu+\frac{1}{2}}e^{-\frac{1}{2}x^2}L_j^\mu(x^2),
\end{equation*}
which is the well--known formula for the Hankel transform of Laguerre polynomials (see e.g. \cite[8.9~(3)]{EMOT54a}).
\end{reprodformula}

\addcontentsline{toc}{section}{References}
\providecommand{\bysame}{\leavevmode\hbox to3em{\hrulefill}\thinspace}
\providecommand{\MR}{\relax\ifhmode\unskip\space\fi MR }
\providecommand{\MRhref}[2]{\href{http://www.ams.org/mathscinet-getitem?mr=#1}{#2}}
\providecommand{\href}[2]{#2}

\vspace{30pt}

\textsc{Joachim Hilgert\\Institut f\"ur Mathematik, Universit\"at Paderborn, Warburger Str. 100, 33098 Paderborn, Germany.}\\
\textit{E-mail address:} \texttt{hilgert@math.uni-paderborn.de}\\

\textsc{Toshiyuki Kobayashi}\\
\textit{Home address:} \textsc{Graduate School of Mathematical Sciences, IPMU, the University of Tokyo, 3-8-1 Komaba, Meguro, Tokyo, 153-8914, Japan.}\\
\textit{Current address:} \textsc{Max-Planck-Institut f\"ur Mathematik, Vivatsgasse 7, 53111 Bonn, Germany.}\\
\textit{E-mail address:} \texttt{toshi@ms.u-tokyo.ac.jp}\\

\textsc{Gen Mano\\Graduate School of Mathematical Sciences, the University of Tokyo, 3-8-1 Komaba, Meguro, Tokyo, 153-8914, Japan.}\\

\textsc{Jan M\"ollers\\Institut f\"ur Mathematik, Universit\"at Paderborn, Warburger Str. 100, 33098 Paderborn, Germany.}\\
\textit{E-mail address:} \texttt{moellers@math.uni-paderborn.de}


\begin{thebibliography}{1}

\bibitem{AAR99}
G.~E. Andrews, R.~Askey, and R.~Roy, \emph{Special functions}, Encyclopedia of
  Mathematics and its Applications, vol.~71, Cambridge University Press,
  Cambridge, 1999.

\bibitem{CL55}
E.~A. Coddington and N.~Levinson, \emph{Theory of ordinary differential
  equations}, McGraw-Hill Book Company, Inc., New York, 1955.

\bibitem{EMOT81}
A.~Erd{\'e}lyi, W.~Magnus, F.~Oberhettinger, and F.~G. Tricomi, \emph{Higher
  transcendental functions. {V}ols. {I}, {II}}, McGraw-Hill Book Company, Inc.,
  New York, 1953.

\bibitem{EMOT54a}
\bysame, \emph{Tables of integral transforms. {V}ol. {II}}, McGraw-Hill Book
  Company, Inc., New York, 1954.

\bibitem{HKMM09a}
J.~Hilgert, T.~Kobayashi, G.~Mano, and J.~M{\" o}llers, \emph{Special functions
  associated to a certain fourth order differential equation}, 
to appear in Ramanujan J. Math. (\href{http://arxiv.org/abs/0907.2608}{arXiv:0907.2608})

\bibitem{KM07bprocj}
T.~Kobayashi and G.~Mano,
  \emph{\href{http://www.springerlink.com/content/pt65717043175035/fulltext.pd%
f}{ Integral formula of the unitary inversion operator for the minimal
  representation of {$O(p,q)$}}}, Proc. Japan Acad. Ser. A (2007),
27--31.
\bibitem{KM07b}
 \emph{The Schr\"odinger model
for the minimal representation of the indefinite
orthogonal group $O(p,q)$},
the Mem. Amer. Math. Soc. (2011)
vol. 212, no. 1000 
\href{http://dx.doi.org/10.1090/S0065-9266-2011-00592-7}{DOI: 10.1090/S0065-9266-2011-00592-7}
(available at
\href{http://arxiv.org/abs/0712.1769}{  arXiv:0712.1769}).

\bibitem{Kra40}
H.~L. Krall, \emph{On orthogonal polynomials satisfying a certain fourth order
  differential equation}, Pennsylvania State College Studies (1940), no.~6.

\bibitem{LK89}
L.~L. Littlejohn and A.~M. Krall, \emph{Orthogonal polynomials and higher
  order singular {S}turm-{L}iouville systems}, Acta Appl. Math. \textbf{17}
  (1989), no.~2, 99--170.

\bibitem{Wat44}
G.~N. Watson, \emph{A {T}reatise on the {T}heory of {B}essel {F}unctions},
  Cambridge University Press, Cambridge, England, 1944.

\end{thebibliography}
\end{document}